\documentclass[graybox]{svmult}
\usepackage{mathptmx}       
\usepackage{helvet}         
\usepackage{courier}        

\usepackage{amsmath}
\usepackage{amssymb} 
\usepackage{amscd}
\usepackage{amsfonts}
\usepackage{makeidx}         
\usepackage{graphicx}        

\usepackage{multicol}        
\usepackage[bottom]{footmisc}
\makeindex

\spdefaulttheorem{q}{Question}{\bfseries}{\itshape}
\spdefaulttheorem{conj}{Conjecture}{\bfseries}{\itshape}

\font\eu=eufm10
\def\ss{\hbox{\eu \char 83}}

\def\moins{\mathrel{\hbox{\vrule height 3pt depth -2pt width 6pt}}}
\def\rond{\kern 1pt{\scriptstyle\circ}\kern 1pt}

\newcommand\Pic{\operatorname{Pic}}

\newcommand\rk{\operatorname{rk}}

\def\Z{\mathbb{Z}}

\def\C{\mathbb{C}}
\def\P{\mathbb{P}}
\def\G{\mathbb{G}}
\def\operp{\,\buildrel {\scriptscriptstyle\perp}\over {\oplus}\,}

\def\qfl#1{\buildrel {#1}\over {\longrightarrow}}

\def\iso{\vbox{\hbox to .8cm{\hfill{$\scriptstyle\sim$}\hfill}
\nointerlineskip\hbox to .8cm{{\hfill$\longrightarrow $\hfill}} }}
\def\diso{\vbox{\hbox to .8cm{\hfill{$\scriptstyle\sim$}\hfill}
\nointerlineskip\hbox to .8cm{{\hfill$\dasharrow $\hfill}} }}

\def\operp{\,\buildrel {\scriptscriptstyle\perp}\over {\oplus}\,}

\font\san=cmssdc10
\def\ext{\hbox{\san \char3}}

\begin{document}
\title*{Holomorphic symplectic geometry:\\ a problem list}
\author{Arnaud Beauville}
\institute{Arnaud Beauville \at Laboratoire J.-A. Dieudonn\'e (UMR 6621 du CNRS),
Universit\'e de Nice,
Parc Valrose,\\
F-06108 Nice cedex 2, France\\
\email{arnaud.beauville@unice.fr}}

\maketitle 
\abstract{ The usual structures of symplectic geometry (symplectic, contact, Poisson) make sense for complex manifolds; they turn out to be quite interesting on projective, or compact K\"ahler, manifolds. In these notes we review some of the recent results on the subject, with emphasis on the open problems and conjectures.}

\section*{Introduction}
Though symplectic geometry is usually done on real manifolds, the main definitions (symplectic or contact structures, Poisson bracket) make perfect sense in the holomorphic setting. What is less obvious is that these structures are indeed quite interesting in this set-up, in particular on \emph{global} objects -- meaning compact, or projective, manifolds. The study of these objects has been much developed in the last 30 years -- an exhaustive survey  would require at least a book. The aim of these notes is much more modest: we would like to give a (very partial) overview of the subject by 
presenting some of the open problems which are currently investigated. 
\par  Most of the paper is devoted to holomorphic symplectic (= hyperk\"ahler) manifolds, a subject which has been blossoming in  recent years. Two short  chapters are devoted to contact and Poisson structures : in the former we discuss the conjectural classification of projective contact manifolds, and in the latter an intriguing conjecture of Bondal on the rank of the Poisson tensor.

\section{Compact hyperk\"ahler manifolds}
\subsection{Basic definitions}
The interest for holomorphic symplectic manifolds comes from the following result, stated by Bogomolov in \cite{Bo1} :
\begin{theorem}[Decomposition theorem]
Let $X$ be a compact, simply-connected K\"ahler manifold with trivial canonical bundle. Then 
$X$ is a product of manifolds of the following two types:
\begin{itemize}
\item projective manifolds $Y$ of dimension  $\geq 3$, with $H^0(Y,\Omega^*_Y)=\C\oplus \C\omega $, where $\omega $ is a generator of $K_Y$;
\item  compact K\"ahler manifolds $Z$ 
with  $H^0(Z,\Omega^*_Z)=\C[\sigma ]$, where $\sigma \in H^0(Z,\Omega^2_Z) $ is everywhere non-degenerate. 
\end{itemize}
\end{theorem}
\par This theorem has an important interpretation (and a proof)  in terms of Riemannian geometry\footnote{See \cite{B4} for a more detailed exposition.}. By the fundamental theorem of Yau \cite{Y}, a $n$-dimensional compact K\"ahler manifold $X$ with trivial canonical bundle admits a K\"ahler metric with holonomy group contained in $ \mathrm{SU}(n)$ (this is equivalent to the vanishing of the Ricci curvature). By the Berger and de Rham theorems, $X$ is a product of manifolds with holonomy $\mathrm{SU}(m)$ or $\mathrm{Sp}(r)$; this corresponds to the first and second case of the decomposition theorem. 
\par We will call the manifolds of the first type  \emph{Calabi-Yau} manifolds, and those of the second type 
 \emph{hyperk\"ahler} manifolds (they are also known as \emph{irreducible holomorphic symplectic}). 

\subsection{Examples}
 For Calabi-Yau manifolds we know a huge quantity of examples (in dimension 3, the number of known families approaches 10 000), but relatively 
little general theory. In contrast, we have much information on hyperk\"ahler manifolds, their period map, their cohomology (see below); what is lacking severely is examples.
 In fact, at this time we know two families in each dimension \cite{B1}, and two isolated families in dimension 6 and 10 \cite{OG1}, \cite{OG2} :

\par \smallskip
$a)$
 Let $S$ be a K3 surface. The symmetric product $S^{(r)}:=S^r/\ss_r$   parametrizes subsets of $r$ points in $S$,  counted with multiplicities; it is smooth on the open subset $S^{(r)}_\mathrm{o}$ consisting of subsets with $r$ distinct points, 
but singular otherwise. If we replace ``subset" by (analytic) ``subspace", we obtain a smooth compact manifold, the    
{\it Hilbert scheme} $S^{[r]}$; the natural map 
 $S^{[r]}\rightarrow S^{(r)}$  is an isomorphism above $S^{(r)}_\mathrm{o}$, but it resolves the singularities of $S^{(r)}$. 
 \par Let $\omega $ be a non-zero holomorphic 2-form on $S$. The form $\mathrm{pr}_1^*\omega +\ldots +\mathrm{pr}_r^*\omega $ descends to a non-degenerate 2-form on $S^{(r)}_\mathrm{o}$; it is easy to check that this 2-form extends to a symplectic structure on $S^{[r]}$.
\par\smallskip $b)$ Let $T$ be a 2-dimensional complex torus. The Hilbert scheme $T^{[r]}$  has the same properties as $S^{[r]}$, but it is not simply connected. This is fixed by considering the  composite map $T^{[r+1]}\rightarrow T^{(r+1)}\qfl{s}T$, where $s(t_1,\ldots ,t_r)=t_1+\ldots +t_r$; the fibre $K_r(T):=s^{-1}(0)$ is a hyperk\"ahler manifold of dimension $2r$ (``generalized Kummer manifold").

\par\smallskip $c)$  Let again $S$ be a K3 surface, and $\mathcal{M}$ the moduli space of stable rank 2 vector bundles on $S$, with  Chern classes $c_1=0$, $c_2=4$. According to Mukai \cite{M}, this space has  a holomorphic symplectic structure. It admits a natural compactification $\overline{\mathcal{M}}$, obtained by adding classes of semi-stable torsion free sheaves; it is singular along the boundary, but O'Grady constructs a desingularization of $\overline{\mathcal{M}}$ which is a new  hyperk\"ahler manifold, of dimension 10. 
 
\par\smallskip $d)$ The analogous construction can be done starting from rank 2 bundles  with   $c_1=0$, $c_2=2$ on
 a 2-dimensional complex torus, and taking again some fibre to ensure the simple connectedness. The upshot is a new hyperk\"ahler manifold of dimension 6.
\medskip

 \par In the  two last examples  it would seem simpler to start with a moduli space $\mathcal{M}$ for which the natural compactification $\overline{\mathcal{M}}$ is smooth; in that case $\overline{\mathcal{M}}$  is a  hyperk\"ahler manifold \cite{M}, but it turns out that it is a deformation of $S^{[r]}$ or $K_r(T)$ (G\"ottsche-Huybrechts, O'Grady, Yoshioka ...).  On the other hand, when $\overline{\mathcal{M}}$ is singular, it
 admits a hyperk\"ahler desingularization only in the two cases considered by O'Grady \cite{KLS}. 
 \par Thus it seems that a new idea is required to answer our first problem:
    \begin{q}
Find new examples of hyperk\"ahler manifolds.
\end{q}

\subsection{The period map}
In dimension 2 the only  hyperk\"ahler manifolds are  K3 surfaces;  we know them very well thanks to the \emph{period map}, which associates to a K3 surface $S$ the Hodge decomposition  
\[ H^2(S,\C)=H^{2,0}\oplus H^{1,1}\oplus H^{0,2} \ .\]
This is determined by the position of the line $H^{2,0}$ in
$H^2(S,\C)$ : indeed we have $H^{0,2}=\overline{H^{2,0}}$, and $H^{1,1}$ is the
orthogonal of  $H^{2,0}\oplus  H^{0,2}$ with respect to the intersection form. 
Note that any non-zero element $\sigma $ of $H^{2,0}$ (that is, the class of a non-zero holomorphic 2-form) satisfies $\sigma^2=0$ and $\sigma \cdot \bar\sigma >0$.

\par  To compare the Hodge structures of different K3 surfaces, we consider \emph{marked} surfaces $(S,\lambda )$, where $\lambda $ is an isometry of 
$H^2(S,\Z)$ onto a fixed lattice $L$, the unique even unimodular  lattice $L$ of signature $(3,19)$. 
Then the data of the Hodge structure on $H^2(S,\Z)$ is equivalent to that of the \emph{period point} $\wp(S,\lambda ):=\lambda_{\C}(H^{2,0} )\in \P(L_{\C})$.  By the above remark this point lies in the domain $\Omega \subset \P(L_{\C})$ defined by the conditions $x^2=0$, $x\cdot \bar x>0$.   There is a moduli space $\mathcal{M}_L$ for marked K3 surfaces, which is a \emph{non-Hausdorff} complex manifold;  the \emph{period map} $\wp:\mathcal{M}_L \rightarrow \Omega _L$ is holomorphic. We know a lot about that map, thanks to the work of many people (Piatetski-Shapiro, Shafarevich, Todorov, Siu, ...): 
 \begin{theorem}
 $1)$ \emph{(``local Torelli")} $\wp$ is a local isomorphism.
\par  $2)$ \emph{(``global Torelli")} If $\wp(S,\lambda )=\wp(S',\lambda' )$, $S$ and $S'$ are isomorphic;
\par  $3)$  \emph{(``surjectivity")} Every point of $\Omega $ is the period of some marked K3 surface.
\end{theorem}
\par  Another way of stating 2) is that $S$ and $S'$ are  isomorphic if and only if  
 there is a \emph{Hodge isometry} $H^2(S,\Z)\iso H^2(S',\Z)$ (that is, an isometry inducing an isomorphism of  Hodge structures). There is in fact a more precise statement, see e.g. \cite{BHPV}.
  \par There is a very analogous  picture for  higher-dimensional hyperk\"ahler manifolds. The intersection form is replaced by a canonical quadratic form $q:H^2(X,\Z)\rightarrow \Z$, primitive\footnote{This means that the associated bilinear form  is integral and not divisible by an integer $>1$.}, of signature $(3,b_2-3)$ \cite{B1}. The easiest way to define it is through the Fujiki relation 
\[\int_X\alpha ^{2r}=f_X\,q(\alpha )^r\quad\hbox{for each }\alpha \in H^2(X,\Z)\ ; \]
 this relation determines $f_X$ (the \emph{Fujiki constant} of $X$) and  the form $q$; they depend only on the topological type of $X$.
\par   Let $X$ be a hyperk\"ahler manifold, and $L$ a lattice.  A \emph{marking} of type $L$ of   $X$ is an isometry $\lambda :(H^2(X,\Z),q)\iso L$. The period of $(X,\lambda )$ is the point $\lambda _{\C}(H^{2,0})\in \P(L_{\C})$; as above it belongs to the period domain
\[\Omega _L:=\{[x]\in \P(L_{\C})\ |\ x^2=0\ ,\ x\cdot \bar x>0\}\ . \]
\par Again we have a non-Hausdorff complex manifold $\mathcal{M}_L$ parametrizing hyper\-k\"ahler manifolds of a given dimension with a marking of type $L$; the period map \allowbreak $\wp:\mathcal{M}_L \rightarrow \Omega _L$ is holomorphic. We have:

\begin{theorem} $1)$ The period map $\wp:\mathcal{M}_L\rightarrow \Omega _L$ is a local isomorphism.
\par  $2)$  The restriction of $\wp$ to any  connected component  of $\mathcal{M}_L$ is surjective.
\end{theorem}
\par 1) is proved in {\cite{B1}, and 2) in {\cite{Hu1}. What is missing is the analogue of the global Torelli theorem. 
It has long been known that it cannot hold in the  form given in Theorem 2; in fact, it follows from the results of \cite{Hu1} that any birational map $X\diso X'$ induces a Hodge isometry 
$H^2(X,\Z)\iso H^2(X',\Z)$. This is not the only obstruction:  Namikawa observed \cite{N} that if $T$ is a 2-dimensional complex torus, and $T^*$ its dual torus, the Kummer manifolds $K_2(T)$ and $K_2(T^*)$ (1.2.$b$) have the same period (with appropriate markings), but are not bimeromorphic in general. Thus we can only ask:
\begin{q}
Let $X$, $X'$ be two  hyperk\"ahler manifolds of the same dimension. If there is a Hodge isometry $\lambda :H^2(X,\Z)\iso H^2(X',\Z)$, what can we say of $X$ and $X'$? Can we conclude that $X$ and $X'$ are isomorphic by imposing that $\lambda $ preserves some extra structure?
\end{q}

\par A partial answer to these questions  appear in \cite{V3}, in particular for the case of example 1.2.$a)$.

\subsection{Cohomology}
Let $X$ be a hyperk\"ahler manifold. Since the quadratic form $q$  plays such an important role, it is natural to expect that it determines most of the cohomology of $X$. This was indeed shown by Bogomolov \cite{Bo3} :

\begin{proposition}
Let $X$ be a hyperk\"ahler manifold, of dimension $2r$, and let $\mathcal{H}$ be the subalgebra of $H^*(X,\C)$ spanned by $H^2(X,\C)$.
\par $1)$ $\mathcal{H}$ is the quotient of $\mathrm{Sym}^*H^2(X,\C)$ by the ideal spanned by the classes $\alpha ^{r+1}$ for $\alpha \in H^2(X,\C)$, $q_{\C}(\alpha )=0$.
\par $2)$ $H^*(X,\C)=\mathcal{H}\oplus \mathcal{H}^\perp$, where $\mathcal{H}^\perp$ is the orthogonal of $\mathcal{H}$ with respect to the cup-product.
\end{proposition}
\par  Thus the subalgebra $\mathcal{H}$ is completely determined by the form $q$ and the dimension of $X$. In contrast, not much is known about the $\mathcal{H}$-module $\mathcal{H}^\perp$. Note that it is nonzero for the  examples  $a)$ and $b)$ of 1.2, with the exception of $S^{[2]}$ for a K3 surface $S$.
\par We do not know much about the quadratic form $q$  either. For the two infinite series of (1.2) we have lattice isomorphisms \cite{B1}
\[H^2(S^{[r]},\Z)=H^2(S,\Z)\operp \langle 2-2r\rangle\qquad H^2(K_r(T),\Z)=H^2(T,\Z)\operp \langle -2-2r\rangle\ ;
\]
The lattices of O'Grady's two examples are computed in \cite{R}; they are also even. 
\begin{q}
 Is the quadratic form $q$ always even? More generally, what are the possibilities for  $q$? What are the possibilities for the Fujiki index $f_X$ (see $1.3$)?
\end{q}

\subsection{Boundedness} Having so few examples leads naturally to the following question:

\begin{conj}
There are finitely many hyperk\"ahler manifolds (up to deformation) in each dimension.
\end{conj}
\par  Note that the same question can be asked  for Calabi-Yau manifolds, but there it seems completely out of reach. 
 \par Huybrechts observes that there are finitely many deformation types of hyperk\"ahler manifolds $X$ of  dimension $2r$ such that there exists $\alpha \in H^2(X,\Z)$ with $q(\alpha )>0$ and $\int_X\alpha ^{2r}$ bounded \cite{Hu2}. As a corollary, given a real number $M$, there are finitely many deformation types of hyperk\"ahler manifolds with
 \[f_X\leq M\quad , \quad  \min\{q(\alpha )\ |\ \alpha \in H^2(X,\Z)\ ,\ q(\alpha )>0\}\leq M\ .
 \]
 
 \par A first approximation to finiteness  would be to bound the Betti numbers $b_i$ of $X$, and in particular $b_2$.
 Here we have some more information in the case of 
 fourfolds \cite{G} :
 \begin{proposition}
Let $X$ be a hyperk\"ahler fourfold. Then either $b_2=23$, or $3\leq b_2\leq 8$.
\end{proposition}
\par Note that $b_2$ is $23$ for $S^{[2]}$ and $7$ for  $K_2(T)$ (1.2).
 \cite{G} contains some more information on the other Betti numbers. 
 \begin{q}
Can we exclude some more cases, in particular $b_2=3$? If $b_2=23$, can we conclude that $X$ is deformation equivalent to $S^{[2]}$?
\end{q}
 \subsection{Lagrangian fibrations}
 Let $(X,\sigma) $ be  a holomorphic symplectic manifold (not necessarily compact), of dimension $2r$. A \emph{Lagrangian fibration} is a proper map $h:X\rightarrow B$ onto a manifold $B$ such that the general fibre $F$ of $h$ is Lagrangian, that is, 
$F$ is connected, of dimension  $r$, and $\sigma _{\,|F}=0$. This implies that the  smooth fibres of $h$ are complex tori
(Arnold-Liouville theorem). 
\par Suppose $B=\C^r$, so that  $h=(h_1,\ldots ,h_r)$. The functions $h_i$ define what is called in classical mechanics an \emph{algebraically completely integrable hamiltonian system} : the Poisson brackets $\{h_i,h_j\}$ vanish,  the hamiltonian vector fields $X_{h_i}$ commute with each other, they are tangent to the fibres of $h$ and their restriction to a smooth fibre is a linear vector field on this complex torus (see for instance \cite{B3}).
 \par  The analogue of this notion when $X$ is compact (hence  hyperk\"ahler)  is a  Lagrangian fibration $X\rightarrow \P^r$. There are many examples of such fibrations (see a sample below); moreover they turn out to be the only non-trivial morphisms 
 from a hyperk\"ahler manifold to a manifold of smaller dimension :
\begin{theorem}
Let $X$ be a hyperk\"ahler manifold, of dimension $2r$, $B$ a  K\"ahler manifold  with $0<\dim B<2r$, and  $f:X\rightarrow B$  a surjective morphism with connected fibres. Then:
\par $1)$ $f$ is a Lagrangian fibration;
\par $2)$ If $X$ is projective, $B\cong \P^r$.
\end{theorem}
\par 1) is due to Matsushita (see \cite{Hu3}, Prop. 24.8), and 2) to Hwang \cite{Hw}. It is expected that 2) holds without the projectivity assumption on $X$ (see the discussion in the introduction of \cite{Hw}).

\smallskip

\par  How do we detect the existence of a Lagrangian fibration on a given  hyperk\"ahler manifold? In dimension 2 there is a simple answer;  a Lagrangian fibration on a K3 surface $S$ is an elliptic fibration, and we have  :
\begin{proposition}
{\rm a)} Let $L$ be a  nontrivial nef line bundle on $S$ with $L^2=0$. There exists an elliptic fibration $f:S\rightarrow \P^1$ such that $L=f^*\mathcal{O}_{\P^1}(k)$ for some $k\geq 1$.\label{K3}
\par  {\rm b)}  $S$ admits an elliptic fibration if and only if it admits a  line bundle $L\neq \mathcal{O}_S$  with $L^2=0$.
\end{proposition}
\par  The proof of a) is straightforward.  b) is reduced to a) by proving that some isometry $w$ of $\Pic(S)$ maps $L$ to a nef line bundle; see for instance \cite{BHPV}, VIII, Lemma 17.4.
\par Proposition \ref{K3}  has a natural (conjectural) generalization to higher-dimensional hyperk\"ahler manifolds\footnote{The conjecture has been known to experts for a long time; see the introduction of \cite{V2} for a discussion of its history.} :

\begin{conj}
{\rm a)} Let $L$ be a  nontrivial nef line bundle on $X$ with $q(L)=0$. There exists a Lagrangian fibration $f:X\rightarrow \P^r$ such that $L=f^*\mathcal{O}_{\P^r}(k)$ for some $k\geq 1$.\label{Lag}
\par  {\rm b)} There exists a hyperk\"ahler manifold $X'$ bimeromorphic to $X$ and a Lagrangian fibration $X'\rightarrow \P^r$ if and only if $X$ admits a  line bundle $L\neq \mathcal{O}_S$ with $q(L)=0$.  
\end{conj}
\par  Note that it is not clear whether one of the statements implies the other.
\par There is some evidence in favor of the conjecture. Let $S$ be a ``general" K3 surface of  genus $g$ -- that is, $\Pic(S)=\Z\,[L]$ with $L^2=2g-2$. Then 
$\Pic(S^{[r]})$ is a rank 2 lattice with an orthonormal basis $(h,e)$ satisfying
$q(h)=2g-2$, $q(e)=-(2r-2)$ \cite{B1}.
 Taking $r=g$ we find $q(h\pm e)=0$. The corresponding Lagrangian fibration  is studied in \cite{B2}: $S^{[g]}$ is birational to the relative compactified Jacobian $\mathcal{J}^g\rightarrow |L|$, whose fibre above a curve $C\in |L|$ is the compactified Jacobian $\overline{J^gC}$. $\mathcal{J}^g$ is hyperk\"ahler by \cite{M}, and the fibration $\mathcal{J}^g\rightarrow |L|$ is Lagrangian. The rational map 
$S^{[g]}\dasharrow |L|$ associates to a general set of $g$ points in $S$ the unique curve of $|L|$ passing through these points.
\par More generally, suppose  that $2g-2=(2r-2)m^2$ for some integer $m$. Then $q(h\pm me)=0$, and indeed $S^{[r]}$ admits a birational  model with a  Lagrangian fibration. This fibration has been constructed independently in \cite{Mar} and \cite{S}; $\mathcal{J}^g$ is replaced by a moduli space of \emph{twisted sheaves} on $S$.
\smallskip
\par  Another argument in favor of the conjecture has been given  by Matsushita \cite{Ma}, who  proved that  b)  holds ``locally", in the following sense. 
 Let $X$ be a hyperk\"ahler manifold, with a Lagrangian fibration $f:X\rightarrow \P^r$, and let $\mathrm{Def}(X)$ be the local deformation space of $X$. Then the Lagrangian fibration deforms along a hypersurface in $\mathrm{Def}(X)$. Thus any small deformation of $X$ such that the cohomology class of  $f^*\mathcal{O}_{\P^r}(1)$  remains algebraic carries a Lagrangian fibration.
\smallskip
\par A related question, which comes from mathematical physics, is :
\begin{q}
Does every hyperk\"ahler manifold admit a deformation with a Lagrangian fibration?
\end{q}
\par If Conjecture \ref{Lag} holds,  the answer is positive if and only if the quadratic form $q$ is indefinite. I do not know any serious argument either in favor or against this. 
\smallskip
\par  
\begin{q}
Let $X$ be a hyperk\"ahler manifold, and $T\subset X$  a Lagrangian submanifold which is a complex torus. Is it the fibre of a  Lagrangian fibration $X\rightarrow \P^r$? 
\end{q}
 \par   (A less optimistic version would ask only for a bimeromorphic Lagrangian fibration.)
 \subsection{Projective families}

 \par Deformation theory shows that when the K3 surface  $S$ varies, the manifolds $S^{[r]}$  form a hypersurface in their deformation space; thus 
  a general deformation of $S^{[r]}$ is not the Hilbert scheme of a K3 -- and 
  we do not know how to describe it. This is not particularly surprising: after all, we do not know either how to describe a general K3 surface. On the other hand, if we start from the family of \emph{polarized} K3 surfaces $S$ of genus $g$, 
 the projective manifolds $S^{[r]}$ are  polarized (in various ways)\footnote{For $S$ general we have $\Pic(S^{[r]})=\Z h\operp \Z e$ (1.6);  the polarizations on $S^{[r]}$ are of the form $ah-be$ with $a,b>0$.}, and the  same argument tells us that they
  form again a hypersurface in their (polarized) deformation space; we should be able to describe 
    a (locally) complete family of projective hyperk\"ahler manifolds which specializes to $S^{[r]}$ in codimension 1. 
    
     \par     For $r=2$ there are indeed a few cases where we can describe the general  deformation of $S^{[2]}$ with an appropriate polarization : 
 \begin{enumerate}

 \item  The Fano variety of lines contained in a cubic fourfold (\cite{BD}; $g=8$)
 
 \item  The ``variety of sum of powers" associated to a cubic fourfold (\cite{IR1}; $g=8$)
  
   \item The double cover of certain sextic hypersurfaces in $\P^5$ (\cite{OG3}; $g=6$)
  
 \item 
  The subspace of the Grassmannian $\G(6,V)$ consisting of $6$-planes $L$ such that $\sigma_{|\,L}=0 $, where $\sigma:\wedge^3\C^{10} \rightarrow \C$ is a sufficiently general 3-form (\cite{DV}; $g=12$).
   \end{enumerate}
\par  Note that K3 surfaces of genus 8 appear  in both cases 1) and 2); what happens is that the corresponding polarizations  on $S^{[2]}$ are different  \cite{IR2}\footnote{The Corollary in \cite{IR2} is slightly misleading: the moduli spaces of \emph{polarized} hyperk\"ahler manifolds of type 1) and of type 2) are disjoint.}.

\begin{q}
Describe the general projective deformation of $S^{[2]}$, for $S$ a polarized K3 surface of genus $1$, $2$, $3$, ...  (and for some choice of polarization on $S^{[2]}$); or at least find more examples of locally complete projective families. Same question with $S^{[r]}$ for $r\geq 3$.
\end{q} 
 \par  (With the notation of footnote 4, a natural choice of polarization  for $g\geq 3$ is $h-e$.)
 
 \medskip
 \par A  different issue concerns the \emph{Chow ring} of a projective hyperk\"ahler manifold. In \cite{B5} and \cite{Vo} the following conjecture is proposed :
 \begin{conj}
Let $D_1,\ldots ,D_k$ in $\Pic(X)$, and let $z\in CH(X)$ be a class which is a polynomial in $D_1,\ldots ,D_k$ and the Chern classes $c_i(X)$. If $z=0$ in $H^*(X,\Z)$, then $z=0$.\label{chow}
\end{conj}
\par This would follow from a much more general (and completely out of reach) conjecture, for which we refer to the introduction of \cite{B5}. Conjecture \ref{chow} is proved in \cite{Vo} for the Hilbert scheme $S^{[n]}$ of a K3 surface for $n\leq 8 $, and for the Fano variety of lines on a cubic fourfold.

\font\cmex=cmex10
\def\ex{\hbox{\cmex \char 96}}
\def\sd_#1{\raise9pt\hbox{$\ \mathop{\kern0pt\ex}\limits_{#1}\ $}}

\section{Compact Poisson manifolds}
 \par Since hyperk\"ahler manifolds are so rare, it is natural to turn to a more flexible notion. Symplectic geometry provides a natural candidate, {\it Poisson manifolds}. Recall that a (holomorphic) Poisson structure on a complex manifold $X$ is a bivector field $\tau \in H^0(X,\ext^2T_X)$, such that the bracket $\{f,g\}:=\langle \tau ,df\wedge dg\rangle$ defines a Lie algebra structure on $\mathcal{O}_X$. A Poisson structure defines a skew-symmetric  map $\tau ^\sharp:\Omega ^1_X\rightarrow T_X$; the {\it rank} of $\tau $ at a point $x\in X$ is the rank of $\tau ^\sharp(x)$. It is  even (because $\tau ^\sharp$ is skew-symmetric). The data of a Poisson structure of rank $\dim X$ is equivalent to that of a (holomorphic) symplectic structure. In general, we have a partition
 
 \[X=\sd_{s\ \mathrm{ even} }X_s \quad\hbox{where }\ X_s :=\{x\in X\ |\ \rk \tau (x)=s\}\ .
 \]

 \par The following conjecture is due to Bondal (\cite{Bd}, see also \cite{P}):
 \begin{conj}
 If $X$ is Fano and $s$ even, $X_{\leq s}:=\sd_{k\leq s}X_k$  contains a component of dimension $>s$. 
\end{conj}
\par This is much larger than one would expect from a naive dimension count. It implies for instance that  a Poisson field which vanishes at some point  must vanish along a curve. 
\par The condition ``$X$ Fano" is probably far too strong. In fact an optimistic modification would be :
 \begin{conj}
 If $X_s$ is non-empty, it  contains a component of dimension $>s$. 
\end{conj}

\par Here are some arguments in favor of this conjecture:
\begin{proposition} Let $(X,\tau )$ be a compact Poisson manifold.
\par  $1)$ Every component  of $X_s$ has dimension $\geq s$.
\par  $2)$
Let $r$ be the generic rank of $\tau $ ($r$ even); assume that $c_1(X)^q\neq 0$ in $H^q(X,\Omega ^q_X)$, where $q=\dim X-r+1$.
Then the  degeneracy locus $X\moins X_{r}$ of $\tau$  has a component of dimension $>r-2$. 
\par  $3)$ Assume that $X$ is a projective threefold. If $X_0$  is non-empty, it contains a curve.
\end{proposition}
\textit{Sketch of proof} : 1) Let $Z$ be a component of $X_s$ (with its reduced structure). It is not difficult to prove that $Z$ is a Poisson subvariety of $X$ (see \cite{P}); this means that at a smooth point $x$ of $Z$, the tensor $\tau (x)$ lives in $\ext^2T_{x}(Z)\subset \ext^2T_{x}(X)$. But this implies $s\leq \dim T_{x}(Z)=\dim Z$. 
\par  2) is proved in  \cite{P}, \S 9, under the extra hypothesis $\dim X= r+1$. The proof extends easily to the slightly more general situation considered here.
\par  3) is proved in \cite{D} by a case-by-case analysis (leading to a complete classification of those Poisson threefolds for which $X_0=\varnothing $).
It would be interesting to have a more conceptual proof.\smallskip
\par The paper \cite{P} contains many interesting results on Poisson manifolds; in particular, a complete classification of the Poisson structures on $\P^3$ for which the zero locus contains a smooth curve.

\section{Compact contact manifolds}
 Let $X$ be a complex manifold, of odd dimension $2r+1$. A {\it contact structure} on $X$ is a one-form $\theta $ with values in a line bundle $L$ on $X$, such that $\theta \wedge (d\theta)^r \neq 0$ at each point of $X$ (though $\theta $ is a twisted 1-form, it is easy to check that $\theta \wedge (d\theta)^r$ makes sense as a section  of  $K_X\otimes L^{r+1}$; in particular, the condition on $\theta $ implies $K_X=L ^{-r-1}$). 

 \par There are only two classes of compact holomorphic contact manifolds known so far:
 \par a) The projective cotangent bundle $\P T^*_M$, where $M$ is any compact complex manifold;
 \par b) Let $\mathfrak{g}$ be a simple complex Lie algebra. The action of the adjoint group on $\P (\mathfrak{g})$ has a unique closed orbit $X_\mathfrak{g}$: every other orbit contains $X_\mathfrak{g}$ in its closure.    $X_\mathfrak{g}$ is a contact Fano manifold.
 \par The following conjecture is folklore:
 \begin{conj}
Any projective contact manifold is of  type {\rm a)} or  {\rm b)}.\label{contact}
\end{conj}
\par Half of this conjecture is now proved, thanks to \cite{KPSW} and \cite{D}:  a contact projective manifold  is either Fano with $b_2=1$, or of type a). It is easily seen that a homogeneous Fano contact manifold is of type b), so  we can rephrase Conjecture \ref{contact} as :
 \begin{conj}
 A contact Fano manifold is homogeneous.\label{Fano}
\end{conj}
\par 
I refer to \cite{B3} for some evidence in favor of this conjecture, and to \cite{B4} for its application to differential geometry, more specifically to   quaternion-K\"ahler manifolds. These are
 Riemannian manifolds with holonomy $\mathrm{Sp}(1)\mathrm{Sp}(r)$;  they are  Einstein mani\-folds, and in particular they have constant scalar curvature. Thanks to work of Salamon  and LeBrun \cite{L,LS}, a positive answer to Conjecture \ref{Fano} would imply:
 \begin{conj}
The only compact quaternion-K\"ahler manifolds with positive scalar curvature are homogeneous.\label{qk}
\end{conj}
 \smallskip
 \par These positive homogeneous quaternion-K\"ahler manifolds have been classified by Wolf \cite{W} : there is one, $M_\mathfrak{g}$, for each simple complex Lie algebra $\mathfrak{g}$. 
 \par The link between Conjectures  \ref{Fano} and \ref{qk} is provided by the twistor space construction. 
To any  quaternion-K\"ahler manifold $M$ is associated
 a $\mathbb{S}^2$-bundle $X\rightarrow M$, the  twistor space, which carries a natural complex structure; when $M$ has positive scalar curvature it turns out that $X$ is a contact Fano manifold -- for instance the twistor space of $M_\mathfrak{g}$ is $X_\mathfrak{g}$. 
Conjecture \ref{Fano} implies that $X$ is isomorphic to $X_\mathfrak{g}$ for some simple Lie algebra $\mathfrak{g}$;  this in turn implies that $M$ is isometric to $M_\mathfrak{g}$ and therefore homogeneous.

\end{document}